\documentclass[twoside,11pt,reqno]{article}
\usepackage{amssymb}
\usepackage{amsthm}
\usepackage[tbtags]{amsmath}
\usepackage{doc}
\usepackage{latexsym}
\usepackage{amscd}
  \theoremstyle{change}

\theoremstyle{change}  
\begin{document}
\thispagestyle{plain}
 \markboth{}{}
\small{\addtocounter{page}{0} \pagestyle{plain}
 \vspace{0.2in}
\noindent{\large \bf Carlitz $q$-Bernoulli Numbers and $q$-Stirling
Numbers}
\vspace{0.15in}\\
\noindent{\sc Taekyun Kim}
\newline
{\it EECS, Kyungpook National University, Daegu 702-701, Korea\\
e-mail} : {\verb|tkim@knu.ac.kr|}
\vspace{0.15in}\\
{\footnotesize {\sc Abstract.} In this paper, we consider Carlitz
$q$-Bernoulli numbers and $q$-stirling numbers of the first and the
second kind. From the properties of $q$-stirling numbers, we derive
many interesting formulae associated with Carlitz $q$-Bernoulli
numbers. Finally, we will prove
$$\beta_{n, q}=\sum_{m=0}^{n}
\sum_{k=m}^{n} \dfrac{1}{(1-q)^{n+m-k}} \sum_{d_{0}+\cdots
+d_{k}=n-k} q^{\sum_{i=0}^{k}i d_{i}} s_{1, q}(k,
m)(-1)^{n-m}\dfrac{m+1}{[m+1]_{q}},$$ where $\beta_{n, q}$ are
called Carlitz $q$-Bernoulli numbers.}
\vspace{0.2in}\\
\noindent{\bf 1. Introduction} \setcounter{equation}{0}
\vspace{0.1in}\\
\indent Let $p$ be a fixed prime number. Throughout this paper,
$\mathbb{Z}_{p}$, $\mathbb{Q}_{p}$, $\mathbb{C}$ and
$\mathbb{C}_{p}$ will, respectively, denote the ring of $p$-adic
rational integers, the field of $p$-adic rational numbers, the
complex number field, and the completion of algebraic closure of
$\mathbb{Q}_{p}$. For $d$ a fixed positive integer with $(p, d)=1$,
let \begin{eqnarray*} X&=&X_{d}=\varprojlim_N
\mathbb{Z}/dp^{N}\mathbb{Z}, \qquad X_{1} = \mathbb{Z}_{p},\\
X^{\ast} &=& \bigcup_{\substack {0<a<dp \\ (a, p)=1}} a+dp \mathbb{Z}_{p},\\
a+dp^{N}\mathbb{Z}_{p}&=&\{x\in X ~|~ x\equiv a ~({\rm
mod}~dp^{N})\}, \end{eqnarray*} where $a \in \mathbb{Z}$ lies in
$0\leq a <dp^{N}$, see [1-21]. The $p$-adic absolute value in
$\mathbb{C}_{p}$ is normalized so that $|p|_{p}=1/p$. When one talks
of $q$-extension, $q$ is variously considered as an indeterminate, a
complex number $q \in \mathbb{C}$ or a $p$-adic number $q\in
\mathbb{C}_{p}$. If $q\in \mathbb{C}_{p}$, then we assume $|q-1|_{p}
<p^{-\frac{1}{p-1}}$, so that $q^{x}=\exp (x \log q)$ for $|x|_p\leq
1$. We use the notation $[x]_{q}=[x:q]=\dfrac{1-q^{x}}{1-q}$. For $f
\in C^{(1)}(\mathbb{Z}_{p})=\{f ~|~ f' \in C(\mathbb{Z}_{p})\}$, let
us start with the expressions
$$\dfrac{1}{[p^{N}]_{q}} \sum_{0\leq j <p^{N}} q^{j}f(j) =
\sum_{0\leq j < p^{N}} f(j)\mu_{q}(j+p^{N}\mathbb{Z}_{p}), \text{
see [6, 8]},$$ representing $q$-analogue of Riemann sums for $f$.
The $p$-adic $q$-integral of a function $f\in
C^{(1)}(\mathbb{Z}_{p})$ is defined by
$$\int_{X} f(x) d\mu_{q}(x) = \int_{\mathbb{Z}_{p}} f(x) d\mu_{q}(x) =
\lim_{N\rightarrow \infty}\dfrac{1}{[p^{N}]_{q}}
\sum^{p^{N}-1}_{x=0} f(x)q^{x}, ~~~~ {\rm see}~[8].$$ For $f \in
C^{(1)}(\mathbb{Z}_{p})$, it is easy to see that,
$$| \int_{\mathbb{Z}_{p}} f(x) d\mu_{q}(x) |_{p} \leq p\|f
\|_{1}, ~~~~ {\rm see}~[6-14],$$ where $\|f \|_{1} = \sup \left\{
|f(0)|_{p}, ~\sup_{x \neq y} |\dfrac{f(x)-f(y)}{x-y}|_{p}\right\}.$
 If $f_{n} \rightarrow f$ in $C^{(1)}(\mathbb{Z}_{p})$, namely $\| f_{n} - f\|_{1}
 \rightarrow 0$, then
 $$\int_{\mathbb{Z}_{p}} f_{n}(x) d\mu_{q}(x) \rightarrow
 \int_{\mathbb{Z}_{p}} f(x) d\mu_{q}(x), ~ {\rm see}~[6-10].$$
 The $q$-analogue of binomial coefficient was known as
 ${\begin{bmatrix}x\\n\end{bmatrix}}_{q} =
 \dfrac{[x]_{q}[x-1]_{q}\cdots[x-n+1]_{q}}{[n]_q !},$ where
 $[n]_{q}!=\prod_{i=1}^{n} [i]_{q},$ (see [1, 5, 6, 10, 11]). From
 this definition, we derive,
 $$\begin{bmatrix}x+1\\{n}\end{bmatrix}_{q} = {\begin{bmatrix}x\\n-1\end{bmatrix}}_{q} + q^{x}
 {\begin{bmatrix}x\\n\end{bmatrix}}_{q} =
 q^{x-n}{\begin{bmatrix}x\\n-1\end{bmatrix}}_{q}  + {\begin{bmatrix}x\\n\end{bmatrix}}_{q} ,~~{\rm cf.~ [6,
 10]}.$$
 Thus, we have
 $\int_{\mathbb{Z}_{p}}{\begin{bmatrix}x\\n\end{bmatrix}}_{q}  d\mu_{q}(x) =
 \dfrac{(-1)^{n}}{[n+1]_{q}} q^{n+1-\binom{n+1}{2}}.$ If
 $f(x)=\underset{k\geq 0}{\sum}a_{k, q}
 {\begin{bmatrix}x\\k\end{bmatrix}}_{q}$ is the $q$-analogue of
 Mahler series of strictly differentiable function $f$, then we see
 that $$\int_{\mathbb{Z}_{p}} f(x) d\mu_{q}(x) = \sum_{k\geq 0} a_{k,
 q}\dfrac{(-1)^{k}}{[k+1]_{q}} q^{k+1-\binom{k+1}{2}}.$$
 Carlitz $q$-Bernoulli numbers $\beta_{k, q}(=\beta_{k}(q))$ can
 be determined inductively by
 $$ \beta_{0, q}=1, ~~q(q\beta +1)^{k} - \beta_{k, q}= \begin{cases}
 1 &{\rm if} ~k=1 \\ 0 &{\rm if} ~k>1, \end{cases}$$ with the usual
 convention of replacing $\beta^{i}$ by $\beta_{i, q}$, (see [2, 3,
 4]). In this paper, we study the $q$-stirling numbers of the first
 and the second kind. From these $q$-stirling numbers, we derive some
 interesting $q$-stirling numbers identities associated with Carlitz
 $q$-Bernoulli numbers. Finally we will prove the following formula
 :
 $$\beta_{n, q} = \sum_{m=q}^{n} \sum_{k=m}^{n}
 \dfrac{1}{(1-q)^{n+m-k}} \sum_{d_{0}+\cdots +d_{k}=n-k}
 q^{\sum_{i=0}^{k} i d_{i}} s_{1, q}(k, m)
 (-1)^{n-m}\dfrac{m+1}{[m+1]_{q}},$$ where $s_{1,q}(k, m)$ is the $q$-stirling
 number of the first kind.
\vspace{0.2in}\\

\noindent{\bf 2. $q$-Stirling numbers and Carlitz $q$-Bernoulli
numbers} \setcounter{equation}{0}
\vspace{0.1in}\\
\indent For $m \in \mathbb{Z}_{+}$, we note that $$\beta_{m,
q}=\int_{\mathbb{Z}_{p}} [x]_{q}^{m} d\mu_{q}(x) = \int_{X}
[x]_{q}^{m} d\mu_{q}(x).$$ From this formula, we derive  $$
\beta_{0, q}=1, ~~q(q\beta +1)^{k} - \beta_{k, q}= \begin{cases}
 1 &{\rm if} ~k=1 \\ 0 &{\rm if} ~k>1, \end{cases}$$ with the usual
 convention of replacing $\beta^{i}$ by $\beta_{i, q}$. By the
 simple calculation of $p$-adic $q$-integral on $\mathbb{Z}_{p}$, we
 see that
 \begin{equation}
 \beta_{n, q}=\dfrac{1}{(1-q)^{n}} \sum_{i=0}^{n} \binom{n}{i}
 (-1)^{i} \dfrac{i+1}{[i+1]_{q}},
 \end{equation}
 where $\binom{n}{i} = \dfrac{n!}{i!(n-i)!}=
 \dfrac{n(n-1)\cdots(n-i+1)}{i!}$. Let $F(t)$ be the generating
 function of Carlitz $q$-Bernoulli numbers. Then we have
 \begin{eqnarray}
 F(t)&=&\sum_{n=0}^{\infty} \beta_{n, q}\dfrac{t^{n}}{n!} =
 \sum_{n=0}^{\infty} \lim_{\rho \rightarrow
 \infty}\dfrac{1}{[p^{\rho}]_{q}} \sum^{p^{\rho}-1}_{x=0} q^{x}
 e^{[x]_{q}t}\\
 &=& \sum_{n=0}^{\infty} \dfrac{1}{(1-q)^{n}}
 \left\{\sum_{k=0}^{\infty} \binom{n}{k}
 \dfrac{k+1}{[k+1]_{q}}(-1)^{k} \right\} \dfrac{t^{n}}{n!}\notag\\
 &=& e^{\dfrac{t}{1-q}} \sum_{k=0}^{\infty}
 \dfrac{(-1)^{k}}{(1-q)^{k}}\dfrac{k+1}{[k+1]_{q}}\dfrac{t^{k}}{k!}.\notag
 \end{eqnarray}
From (2) we note that, \begin{eqnarray}
F(t)&=&e^{\dfrac{t}{1-q}}+e^{\dfrac{t}{1-q}}\sum_{k=1}^{\infty}\dfrac{(-1)^{k}}{(1-q)^{k-1}}
\left(\dfrac{k}{1-q^{k+1}}\right)\dfrac{t^{k}}{k!} \\
~~&&+
e^{\dfrac{t}{1-q}}\sum_{k=1}^{\infty}\dfrac{(-1)^{k}}{(1-q)^{k-1}}\left(\dfrac{1}{1-q^{k+1}}\right)\dfrac{t^{k}}{k!}\notag\\
&=&-t\sum_{n=0}^{\infty}q^{2n}e^{[n]_{q} t} +
(1-q)\sum_{n=0}^{\infty} q^{n}e^{[n]_{q}t}.\notag
\end{eqnarray}
Therefore we obtain the following:
\vspace{0.1in}\\
\noindent {\bf Lemma 1.} {\it Let $F(t)=\sum_{n=0}^{\infty}
\int_{\mathbb{Z}_{p}} [x]_{q}^{n} d\mu_{q}(x) \dfrac{t^{n}}{n!}$.
Then we have $$F(t)=-t\sum_{n=0}^{\infty} q^{2n}e^{[n]_{q}t} +
(1-q)\sum_{n=0}^{\infty} q^{n}e^{[n]_{q}t}.$$}

The $q$-Bernoulli polynomials in the variable $x$ in
$\mathbb{C}_{p}$ with $|x|_{p} \leq 1$ are defined by
\begin{equation}
\beta_{n, q}(x)= \int_{\mathbb{Z}_{p}} [x+t]_{q}^{n} d\mu_{q}(t) =
\int_{X} [x+t]_{q}^{n} d\mu_{q}(x).
\end{equation}
Thus we have
\begin{eqnarray*}
\int_{\mathbb{Z}_{p}} [x+t]_{q}^{n} d\mu_{q}(x) &=&
\sum_{k=0}^{n}\binom{n}{k} [x]_{q}^{n-k} q^{kx}
\int_{\mathbb{Z}_{p}}
[t]_{q}^{k} d\mu_{q}(t)\\
&=& \sum_{k=0}^{n}\binom{n}{k} [x]_{q}^{n-k} q^{kx} \beta_{k, q} =
(q^{x}\beta + [x]_{q})^{n}.
\end{eqnarray*}
From (4) we derive
\begin{equation} \int_{\mathbb{Z}_{p}}
[x+t]_{q}^{n} d\mu_{q}(x) = \beta_{n, q}(x) = \dfrac{1}{(1-q)^{n}}
\sum_{k=0}^{n} \binom{n}{k} (-1)^{k} q^{kx} \dfrac{k+1}{[k+1]_{q}}.
\end{equation}
Let $F(t, x)$ be the generating function of $q$-Bernoulli
polynomials. By (5) we see that \begin{equation} F(t,
x)=\sum_{n=0}^{\infty} \beta_{n, q}(x) \dfrac{t^{n}}{n!} =
e^{\dfrac{t}{1-q}} \sum_{k=0}^{\infty} \dfrac{1}{(1-q)^{k}}
q^{kx}(-1)^{k} \dfrac{k+1}{[k+1]_{q}} \dfrac{t^{k}}{k!}.
\end{equation}
From (6) we note that
\begin{equation}
F(t, x) = -t\sum_{n=0}^{\infty} q^{2n+x} e^{[n+x]_{q}t} + (1-q)
\sum_{n=0}^{\infty} q^{n} e^{[n+x]_{q} t}.
\end{equation}
By (4) and (7), we easily see that
\begin{equation}
[m]_{q}^{k-1} \sum_{i=0}^{m-1} q^{i} \beta_{k,
q^{m}}(\dfrac{x+i}{m}) = \beta_{k, q}(x), \quad m \in \mathbb{N}, k
\in \mathbb{Z}_{+}.
\end{equation}
If we take $x=0$ in (8), then we have
$$[n]_q \beta_{n, q} = \sum_{k=0}^{m}\binom{m}{k} \beta_{k, q^{n}}
[n]_{q}^{k} \sum_{j=0}^{n-1} q^{j(k+1)} [j]_{q}^{n-k}.$$ By (2), (6)
and (7), we see that
\begin{equation}
-\sum_{l=0}^{\infty} q^{2l+n} e^{[n+l]_{q} t} +
\sum_{l=0}^{\infty}q^{2l}e^{[l]_{q} t} = \sum_{m=1}^{\infty} (m
\sum_{l=0}^{n-1} q^{2l}[l]_{q}^{m-1}) \dfrac{t^{m-1}}{m!}.
\end{equation}
Note that $\sum_{l=0}^{\infty} q^{2l+n} e^{[n+l]_{q} t} +
\sum_{l=0}^{n} q^{2 l} e^{[l]_{q} t} = \dfrac{1}{t} (F(t, n)-F(t)).$
Thus, we have
\begin{equation}
\sum_{m=0}^{\infty} (\beta_{m, q}(n) - \beta_{m, q})
\dfrac{t^{m}}{m!} = \sum_{m=0}^{\infty} (m \sum_{l=0}^{n-1} q^{2l}
[l]_{q}^{m-1}) \dfrac{t^{m}}{m!}.
\end{equation}
By comparing the coefficients on both sides in (10), we see that
\begin{equation}
\beta_{m, q}(n) - \beta_{m, q} = m \sum_{l=0}^{n-1} q^{2l}
[l]_{q}^{m-1}. \end{equation} Therefore we obtain the following:
\vspace{0.1in}\\
\noindent {\bf Proposition 2.} {\it For $m, n \in \mathbb{N}$, we
have
$$(q-1)\sum_{l=0}^{n-1} q^{l}[l]_{q}^{m} + \sum_{l=0}^{n-1} q^{l}
[l]_{q}^{m-1} = \dfrac{1}{m} \sum_{l=0}^{m-1} \binom{m}{l}
[n]_{q}^{m-l} q^{nl} \beta_{l, q} + (q^{mn}-1) \beta_{m, q}.$$} Now
we consider the $q$-analogue of Jordan factor as follows:
$$[x]_{k, q} = [x]_{q}[x-1]_{q} \cdots [x-k+1]_{q} =
\dfrac{(1-q^{x})(1-q^{x-1}) \cdots (1-q^{x-k+1})}{(1-q)^{k}}.$$ The
$q$-binomial coefficient is defined by
\begin{equation}
{\begin{bmatrix} n \\ k \end{bmatrix}}_{q} =
\dfrac{[n]_{q}!}{[k]_{q}![n-k]_{q}!} = \dfrac{(1-q^{n})(1-q^{n-1})
\cdots (1-q^{n-k+1})}{(1-q)(1-q^{2})\cdots (1-q^{k})},
\end{equation}
where $[n]_{q}!=[n]_{q} [n-1]_{q} \cdots [2]_{q}[1]_{q}.$ The
$q$-binomial formulas are known as
\begin{equation}
\prod_{i=1}^{n} (a+bq^{i-1}) = \sum_{k=0}^{n} {\begin{bmatrix} n \\
k \end{bmatrix}}_{q} q^{\binom{n}{k}} a^{n-k} b^{k},
\end{equation}
and $$\prod_{i=1}^{n} (1-bq^{i-1})^{-1} = \sum_{k=0}^{n} {\begin{bmatrix} n+k-1 \\
k \end{bmatrix}}_{q} b^{k}.$$ The $q$-Stirling numbers of the first
kind $s_{1, q}(n, k)$ and the second kind $s_{2, q}(n, k)$ are
defined as
\begin{equation}
[x]_{n, q} = q^{-\binom{n}{2}} \sum_{l=0}^{n} s_{1, q}(n, l)
[x]_{q}^{l}, \quad n= 0, 1, 2, \cdots,
\end{equation}
and
\begin{equation}
[x]^{n}_{q} = \sum_{k=0}^{n} q^{\binom{k}{2}} s_{2, q}(n, k) [x]_{k,
q}, \quad n= 0, 1, 2, \cdots, ~{\rm see}~[2, ~3, ~6].
\end{equation}
The values $s_{1, q}(n, 1), \quad n=1, 2, 3, \cdots,$ and $s_{2,
q}(n, 2), \quad n=2, 3, \cdots,$ may be deduced from the following
recurrence relation:
 $$s_{1, q}(n, k)=s_{1, q}(n-1, k-1)-[n-1]_{q}s_{1, q} (n-1, k),
 \quad {\rm see}~[2, ~3, ~6],$$
 for $k=1, 2, \cdots, n$, $n=1, 2, \cdots,$ with initial conditions
 $s_{1, q}(0, 0)=1$, $s_{1, q}(n, k)=0$ if $k > n$. For $k=1$, it
 follows that $$s_{1, q}(n, 1) = -[n-1]_{q} s_{1, q} (n-1, 1),
 \quad n= 2, 3, \cdots,$$ and since $s_{1, q}(1, 1)=1$, we have
 $s_{1, q}(n, 1) = (-1)^{n-1}[n-1]_{q}!, \quad n= 1, 2, 3, \cdots.$
 The recurrence relation for $k=2$ reduce to $s_{1, q}(n, 2) +
 [n-1]_{q} s_{1, q}(n-1, 2) = (-1)^{n-2}[n-2]_{q}!, \quad n= 3, 4,
 \cdots.$ By simple calculation, we easily see that
 \begin{eqnarray*}
 \dfrac{(-1)^{n+1} s_{1, q}(n+1, 2)}{[n]_{q}!} -
 \dfrac{(-1)^{n} s_{1, q}(n, 2)}{[n-1]_{q}!} &=& (-1)^{n+1}
 \dfrac{s_{1, q}(n+1, 2)-[n]_{q} s_{1, q}(n, 2)}{[n]_{q}!}\\
 &=& (-1)^{n+1}\dfrac{(-1)^{n+1}[n-1]_{q}!}{[n]_{q}!} =
 \dfrac{1}{[n]_{q}}, \quad n= 2, 3, 4, \cdots.
 \end{eqnarray*}
Thus we have
$$\dfrac{(-1)^{n}s_{1, q}(n, 2)}{[n-1]_{q}!} = \sum_{k=1}^{n-1}
\dfrac{1}{[k]_{q}}.$$ This is equivalent to $s_{1, q}(n, 2) =
(-1)^{n}[n-1]_{q}! \sum_{k=1}^{n-1}\dfrac{1}{[k]_{q}}.$ It is easy
to see that $$\sum_{m=1}^{n} (-1)^{m+1} q^{\binom{m+1}{2}}
{\begin{bmatrix} n+1 \\ m+1 \end{bmatrix}}_{q} \sum_{k=1}^{m}
\dfrac{1}{[k]_{q}} = \sum_{k=1}^{n} (-1)^{k+1}q^{\binom{k+1}{2}}
\dfrac{{\begin{bmatrix} n \\ k \end{bmatrix}}_{q} }{[k]_{q}}.$$ From
this, we derive
\begin{eqnarray*}
\sum_{k=1}^{n} (-1)^{k+1}q^{\binom{k+1}{2}} \dfrac{1}{[k]_{q}}
\left({\begin{bmatrix} n \\ k \end{bmatrix}}_{q} - {\begin{bmatrix} n-1 \\
k \end{bmatrix}}_{q}\right) &=& \sum_{k=1}^{n} (-1)^{k+1}
q^{\binom{k+1}{2}} \dfrac{1}{[k]_{q}} \left(q^{n-k} {\begin{bmatrix}
n-1 \\ k-1 \end{bmatrix}}_{q}\right)\\
=\dfrac{q^{n}}{[n]_{q}} \sum_{k=1}^{n} (-1)^{k+1} q^{\binom{k}{2}}
{\begin{bmatrix} n \\ k \end{bmatrix}}_{q}  &=&
\dfrac{q^{n}}{[n]_{q}}.
\end{eqnarray*}
Note that $\sum_{k=1}^{n}(-1)^{k+1} q^{\binom{k}{2}}
{\begin{bmatrix} n
\\ k \end{bmatrix}}_{q} = -\sum_{k=0}^{n}(-1)^{k} q^{\binom{k}{2}} {\begin{bmatrix} n \\ k \end{bmatrix}}_{q}
+1 = 1.$ Thus, we have $$\sum_{k=1}^{n}(-1)^{k+1}
q^{\binom{k+1}{2}}\dfrac{ {\begin{bmatrix} n
\\ k \end{bmatrix}}_{q}}{[k]_{q}} = \sum_{k=1}^{n-1}(-1)^{k+1} q^{\binom{k+1}{2}}
\dfrac{{\begin{bmatrix} n-1 \\ k \end{bmatrix}}_{q}}{[k]_{q}} +
\dfrac{q^{n}}{[n]_{q}}.$$

Continuning this process, we see that $$\sum_{k=1}^{n}(-1)^{k+1}
q^{\binom{k+1}{2}}\dfrac{ {\begin{bmatrix} n
\\ k \end{bmatrix}}_{q}}{[k]_{q}} =
\sum_{k=1}^{n}\dfrac{q^{n}}{[k]_{q}}.$$

 The $p$-adic $q$-gamma
function is defined as $\Gamma_{p, q}(n) = (-1)^{n} \prod_{\substack
{ 1\leq j<n \\(j, p) = 1}} [j]_{q}.$ For all  $x \in
\mathbb{Z}_{p}$, we have $\Gamma_{p, q}(x+1) = E_{p, q}(x)
\Gamma_{p, q}(x),$ where $E_{p, q}(x) =  \begin{cases}
 -[x]_{q} &{\rm if} ~|x|_{p}=1 \\ -1 &{\rm if} ~|x|_{p}<1.
 \end{cases}$ Thus, we easily see that
 \begin{equation}
 \log \Gamma_{p, q}(x+1) = \log E_{p, q}(x) + \log \Gamma_{p, q}(x).
 \end{equation}
From the differentiating on both sides in (16), we derive
$$\dfrac{\Gamma_{p, q}^{\prime}(x+1)}{\Gamma_{p, q}(x+1)} = \dfrac{\Gamma_{p, q}^{\prime}(x)}{\Gamma_{p,
q}(x)} + \dfrac{E_{p, q}^{\prime}(x)}{E_{p, q}(x)}.$$ Continuning
this process, we have

 $$\dfrac{\Gamma_{p, q}^{\prime}(x)}{\Gamma_{p,
q}(x)} = \left(\sum_{j=1}^{x-1} \dfrac{q^{j}}{[j]_{q}}\right)
\dfrac{\log q}{q-1} +\dfrac{\Gamma_{p, q}^{\prime}(1)}{\Gamma_{p,
q}(1)}.$$ The classical Euler constant is known as
$\gamma=\dfrac{\Gamma^{\prime}(1)}{\Gamma(1)}$. In [15], Koblitz
defined the $p$-adic $q$-Euler constant as
$$\gamma_{p, q}=-\dfrac{\Gamma_{p, q}^{\prime}(1)}{\Gamma_{p, q}(1)}.$$
Therefore, we obtain the following:
\vspace{0.1in}\\
\noindent {\bf Theorem 3.} {\it For $x \in \mathbb{Z}_{p}$, we have}
$$ \sum_{k=1}^{x-1} (-1)^{k+1} q^{\binom{k+1}{2}} \dfrac{{\begin{bmatrix} x-1 \\ k
\end{bmatrix}}_{q}}{[k]_{q}} = \dfrac{q-1}{\log q} \left( \dfrac{\Gamma_{p, q}^{\prime}(x)}{\Gamma_{p,
q}(x)} - \gamma_{p, q}\right).$$ From (5), (12), (14) and (15), we
derive the following theorem:
\vspace{0.1in}\\
\noindent {\bf Theorem 4.} {\it For $n, k \in \mathbb{Z}_{+}$, we
have
$$\beta_{n, q} = \dfrac{1}{(1-q)^{n}} \sum^{n}_{l=0} \binom{n}{l}
(-1)^{l} \sum^{l}_{k=0} (q-1)^{k}{\begin{bmatrix} l \\ k
\end{bmatrix}}_{q} \sum_{m=0}^{k} s_{1, q} (k, m) \beta_{m, q},$$
where $s_{1, q}(k, m)$ is the $q$-Stirling number of the first
kind.}

By simple calculation, we easily see that
\begin{eqnarray*} q^{nt}&=& ([t]_{q}(q-1)+1)^{n} = \sum_{m=0}^{n}
\binom{n}{m} (-1)^{m} (1-q)^{m} [t]_{q}^{m} =  \sum_{k=0}^{n}
(q-1)^{k} q^{\binom{k}{2}} {\begin{bmatrix} n \\ k
\end{bmatrix}}_{q} [t]_{k, q} \\
&=& \sum_{k=0}^{n} (q-1)^{k} {\begin{bmatrix} n \\ k
\end{bmatrix}}_{q} \sum_{m=0}^{k} s_{1, q}(k, m)[t]_{q}^{m} =
\sum_{m=0}^{n} \left(\sum_{k=m}^{n} (q-1)^{k} {\begin{bmatrix} n \\
k\end{bmatrix}}_{q} s_{1, q}(k, m)\right) [t]_{q}^{m}.
\end{eqnarray*}
Thus we note
\begin{equation}
\int_{\mathbb{Z}_{p}} q^{nt} d\mu_{q}(t) = \sum_{m=0}^{n}
\left(\sum_{k=m}^{n}(q-1)^{k} {\begin{bmatrix} n \\ k
\end{bmatrix}}_{q} s_{1, q}(k, m)\right) \beta_{m, q}.
\end{equation}
From the definition of $p$-adic $q$-integral on $\mathbb{Z}_{p}$, we
also derive
\begin{equation}
\int_{\mathbb{Z}_{p}} q^{nt} d\mu_{q}(t) = \sum_{m=0}^{n}
\binom{n}{m} (q-1)^{m} \beta_{m, q}.
\end{equation}
By comparing the coefficients on the both sides of (17) and (18), we
see that
$$\binom{n}{m}(q-1)^{m} = \sum_{k=m}^{n}(q-1)^{k} {\begin{bmatrix} n \\ k
\end{bmatrix}}_{q} s_{1, q}(k, m).$$
Therefore we obtain the following:
\vspace{0.1in}\\
\noindent {\bf Theorem 5.} {\it For $n \in \mathbb{N}, m \in
\mathbb{Z}_{+}$, we have}
$$\binom{n}{m} = \sum_{k=m}^{n} (q-1)^{-m+k} {\begin{bmatrix} n \\ k
\end{bmatrix}}_{q} s_{1, q}(k, m).$$
From Theorem 5, we can also derive the following interesting formula
for $q$-Bernoulli numbers:
\vspace{0.1in}\\
\noindent {\bf Theorem 6.} {\it For $n \in \mathbb{Z}_{+}$, we have}
$$\beta_{n, q} = \dfrac{1}{(1-q)^{n}} \sum_{m=0}^{n}
\left(\sum_{k=m}^{n}(q-1)^{-m+k} {\begin{bmatrix} n \\ k
\end{bmatrix}}_{q} s_{1, q}(k,
m)\right)(-1)^{m}\dfrac{m+1}{[m+1]_{q}}.$$ From the definition of
$q$-binomial coefficient, we easily derive
\begin{equation}
{\begin{bmatrix} x+1 \\ n
\end{bmatrix}}_{q} = {\begin{bmatrix} x \\ n-1
\end{bmatrix}}_{q} + q^{x}{\begin{bmatrix} x \\ n
\end{bmatrix}}_{q} = q^{x-n} {\begin{bmatrix} x \\ n-1
\end{bmatrix}}_{q} + {\begin{bmatrix} x \\ n
\end{bmatrix}}_{q}.
\end{equation}
By (19), we see that
\begin{equation}
\int_{\mathbb{Z}_{p}} {\begin{bmatrix} x \\ n
\end{bmatrix}}_{q} d\mu_{q}(x) = \dfrac{(-1)^{n}}{[n+1]_{q}}
q^{n+1-\binom{n+1}{2}}.
\end{equation}
From the definition of $q$-Stirling number of the first kind, we
also note that
\begin{equation}
\int_{\mathbb{Z}_{p}} [x]_{n, q}d\mu_q(x)=[n]_{q}!
\int_{\mathbb{Z}_{p}} {\begin{bmatrix} x \\ n
\end{bmatrix}}_{q} d\mu_{q}(x) = q^{-\binom{n}{2}} \sum_{k=0}^{n}
s_{1, q}(n, k)\beta_{k, q}.
\end{equation}
By using (20), (21), we see
\begin{equation}
(-1)^{n} \dfrac{q[n]_{q}!}{[n+1]_{q}} = \sum_{k=0}^{n} s_{1, q}(n,
k) \beta_{k, q}.
\end{equation}
From (15) and (21), we derive
 $$\beta_{n, q} = q \sum_{k=0}^{n}
s_{2, q}(n, k)(-1)^{k} \dfrac{[k]_{q}!}{[k+1]_{q}}.$$ Therefore we
obtain the following:
\vspace{0.1in}\\
\noindent {\bf Theorem 7.} {\it For $ n\in \mathbb{Z}_{+}$, we have
$$ \beta_{n, q} = q \sum_{k=0}^{n}
s_{2, q}(n, k)(-1)^{k} \dfrac{[k]_{q}!}{[k+1]_{q}},$$ where $s_{2,
q}(n, k)$ is the $q$-Stirling number of the second kind.}

It is easy to see that
\begin{equation}
{\begin{bmatrix} n \\ k
\end{bmatrix}}_{q} = \sum_{d_{0}+\cdots + d_{k} = n-k}
q^{\sum_{i=0}^{k} i d_{i}}.
\end{equation}
By Theorem 4, we have the following:
\vspace{0.1in}\\
\noindent {\bf Theorem 8.} {\it For $ n\in \mathbb{Z}_{+}$, we have
$$ \beta_{n, q} = \sum_{m=0}^{n} \sum_{k=m}^{n}
\dfrac{1}{(1-q)^{n+m-k}} \sum_{d_{0}+\cdots+d_{k} = n-k}
q^{\sum_{i=0}^{k} i d_{i}} s_{1, q}(k, m)(-1)^{n-m}
\dfrac{m+1}{[m+1]_{q}},$$ where $s_{1, q}(k, m)$ is the $q$-Stirling
number of the first kind.}
\vspace{0.2in}\\
\footnotesize{

\end{document}